\documentclass[dvips, 10pt, a4paper]{article}
\usepackage{latexsym, amsmath}
\usepackage{amssymb}

\begin{document}

\begin{center}
{\Large{Classification of totally umbilical slant submanifolds of a Kenmotsu manifold}}

\bigskip
\noindent
{\it{Siraj Uddin, Zafar Ahsan and  A.H. Yaakub}}\\
\end{center}

\bigskip
\begin{abstract}
The purpose of this paper is to classify totally umbilical slant submanifolds of a Kenmotsu manifold. We prove that a totally umbilical slant submanifold $M$ of a Kenmotsu manifold $\bar M$ is either invariant or anti-invariant or $dim M=1$ or the mean curvature vector $H$ of $M$ lies in the invariant normal subbundle. Moreover, we find with an example that every totally umbilical proper slant submanifold is totally geodesic.
\end{abstract}
\section{Introduction}

\parindent=8mm
Slant submanifolds of an almost Hermitian manifold were defined by Chen as a natural generalization of both holomorphic and totally real submanifolds [6]. On the other hand, A. Lotta [15] has introduced the notion of slant immersions into almost contact metric manifolds and obtained the results of fundamental importance. He has also studied the intrinsic geometry of $3-$dimensional non anti-invariant slant submanifolds of $K-$contact manifolds [16]. Later on, Cabrerizo et. al [3] studied the geometry of slant submanifolds in more specialized settings of $K-$contact and Sasakian manifolds and obtained many interesting results.

\parindent=8mm On the other hand, in 1954, J.A. Schouten studied the totally umbilical submanifolds and proved that every totally umbilical submanifold of $dim\geq4$ in a conformally flat space is conformally flat [17]. After that many authors studied the geometrical aspects of these submanifolds in different settings, including those of [1, 4, 5, 8, 9, 18]. In this paper, we consider $M$, a totally umbilical slant submanifold tangent to the structure vector field $\xi$ of a Kenmotsu manifold $\bar M$ and obtain a classification result that either $(i)~~M$ is anti-invariant or $(ii)~~dim M=1$ or $(iii)~~H\in\Gamma(\mu)$, where $\mu$ is the invariant normal subbundle under $\phi$. We also prove that every totally umbilical proper slant submanifold is totally geodesic. To, this end, we provide an example to justifiy our results.

\section{Preliminaries}

\parindent=8mm
A $(2n+1)-$dimensional manifold  $(\bar M, g)$ is said to be an {\it{almost contact metric manifold}} if it admits an endomorphism $\phi$ of its tangent bundle $T\bar M$, a vector field $\xi$, called {\it{structure vector field}} and $\eta$, the dual $1-$form of $\xi$ satisfying the following [2]:

$$\phi^2=-I+\eta\otimes\xi,~~\eta(\xi)=1,~~\phi(\xi)=0,~~\eta\circ\phi=0\eqno(2.1)$$
and
$$g(\phi X, \phi Y)=g(X, Y)-\eta(X)\eta(Y),~~\eta(X)=g(X, \xi)\eqno(2.2)$$
for any  $X, Y$ tangent to  $\bar M$. An almost contact metric manifold is known to be {\it{Kenmotsu manifold}} [12] if
$$(\bar\nabla_X\phi)Y=g(\phi X, Y)\xi-\eta(Y)\phi X\eqno(2.3)$$
and
$$\bar\nabla_X\xi=X-\eta(X)\xi\eqno(2.4)$$
for any vector fields $X, Y$ on $\bar M$, where $\bar\nabla$ denotes the Riemannian connection with respect to $g$.

\parindent=8mm
Now, let $M$ be a submanifold of $\bar M$. We will denote by $\nabla$, the induced Riemannian connection on $M$ and $g$, the Riemannian metric on $\bar M$ as well as the metric induced on $M$. Let $TM$ and$T^{\perp}M$ be the Lie algebra of vector fields tangent to $M$ and normal to $M$, respectively and $\nabla^\perp$ the induced connection on $T^\perp M$. Denote by ${\cal{F}}(M)$ the algebra of smooth functions on $M$ and by $\Gamma(TM)$ the ${\cal{F}}(M)$-module of smooth sections of $TM$ over $M$. Then the Gauss and Weingarten formulas are given by
$$\bar\nabla_XY=\nabla_XY+h(X, Y)\eqno(2.5)$$
$$\bar\nabla_XN=-A_NX+\nabla^{\perp}_XN,\eqno(2.6)$$
for each $X,~Y\in\Gamma(TM)$ and $N\in \Gamma(T^\perp M)$, where $h$ and $A_N$ are the second
fundamental form and the shape operator (corresponding to the normal vector
field $N$) respectively for the immersion of $M$ into $\bar M$. They are related as
$$g(h(X, Y), N)=g(A_NX, Y).\eqno(2.7)$$

\parindent=8mm
Now, for any $X\in\Gamma(TM)$, we write
$$\phi X=TX+FX,\eqno(2.8)$$
where $TX$ and $FX$ are the tangential and normal components of $\phi X$, respectively. Similarly for any $N\in\Gamma(T^\perp M)$, we have
$$\phi N=tN+fN,\eqno(2.9)$$
where $tN$ (resp. $fN$) is the tangential (resp. normal) component of $\phi N$.

\parindent=8mm From (2.1) and (2.8), it is easy to observe that for each $X, Y\in\Gamma(TM)$
$$g(TX, Y)=-g(X, TY).\eqno(2.10)$$

\parindent=8mm The covariant derivatives of the endomorphisms $\phi,~T$ and $F$ are defined respectively as
$$(\bar\nabla_X\phi)Y=\bar\nabla_X\phi Y-\phi\bar\nabla_XY,~~\forall X, Y\in\Gamma(T\bar M)\eqno(2.11)$$
$$(\bar\nabla_XT)Y=\nabla_XTY-T\nabla_XY,~~\forall X, Y\in\Gamma(TM)\eqno(2.12)$$
$$(\bar\nabla_XF)Y=\nabla_{X}^{\perp}FY-F\nabla_XY~~\forall X, Y\in\Gamma(TM).\eqno(2.13)$$

\parindent=8mm Throughout, the structure vector field $\xi$ assumed to be tangential to $M$, otherwise $M$ is simply anti-invariant [15]. For any $X\in\Gamma(TM)$, on using (2.4) and (2.5), we may obtain
$$(a)~~\nabla_X\xi=X-\eta(X)\xi,~~~~(b)~~h(X,\xi)=0.\eqno(2.14)$$
On using (2.3), (2.5), (2.6), (2.8), (2.9) and (2.11)-(2.13), we obtain
$$(\bar\nabla_XT)Y=g(TX, Y)\xi-\eta(Y)TX+A_{FY}X+th(X, Y)\eqno(2.15)$$
$$(\bar\nabla_XF)Y=fh(X, Y)-h(X, TY)-\eta(Y)FX.\eqno(2.16)$$

\parindent=8mm A submanifold $M$ of an almost contact metric manifold $\bar M$ is said to be {\it{totally umbilical}} if
$$h(X, Y)=g(X, Y)H,\eqno(2.17)$$
where $H$ is the mean curvature vector of $M$. Furthermore, if $h(X, Y)=0$, for all $X, Y\in\Gamma(TM)$, then $M$ is said to be {\it{totally geodesic}} and if $H=0$, them $M$ is {\it{minimal}} in $\bar M$.

\parindent=8mm For a totally umbilical submanifold $M$ tangent to the structure vector field $\xi$ of a Kenmotsu manifold $\bar M$, we have
$$g(X, \xi)H=0,~~\forall X\in\Gamma(TM).\eqno(2.18)$$
There are two possible cases arise, hence we conclude the following:

\parindent=0mm {\it{Case (i):}} When $X$ and $\xi$ are linearly dependent, i.e., $X=\alpha\xi$, for some non-zero $\alpha\in\mathbb R$, then $g(X, \xi)=\alpha.$ In this case, from (2.18), we get $H=0$ with $dim M=1$, which is trivial case of totally geodesic submanifold of unit dimension.

\parindent=0mm {\it{Case (ii):}} When $X$ and $\xi$ are orthogonal, then from (2.18), it is not necessary that $H=0$, which is the case has to be discussed for totally umbilical submanifolds.

\parindent=8mm In the following section, we will discuss all possible cases of totally umbilical slant submanifolds.

\section{Slant submanifolds}
A submanifold $M$ tangent to the structure vector filed $\xi$ of an almost contact metric manifold $\bar M$ is said to be {\it{slant submanifold}} if for any $x\in M$ and $X\in T_xM-\langle\xi\rangle$, the angle between $\phi X$ and $T_xM$ is constant. The constant angle $\theta\in[0,\pi/2]$ is then called {\it{slant angle}} of $M$ in $\bar M$. Thus, for a slant submanifold $M$, the tangent bundle $TM$ is decomposed as
$$TM=D\oplus\langle\xi\rangle$$
where the orthogonal complementary distribution $D$ of $\langle\xi\rangle$ is known as {\it{slant distribution}} on $M$. The normal bundle $T^\perp M$ of $M$ is decomposed as
$$T^\perp M=F(TM)\oplus\mu,$$
where $\mu$ is the invariant normal subbundle with respect to to $\phi$ orthogonal to $F(TM)$.

\parindent=8mm For a proper slant submanifold $M$ of an almost contact metric manifold $\bar M$ with the slant angle $\theta$, Lotta [15] proved that
$$T^2X=-\cos^2\theta(X-\eta(X)\xi)\eqno(3.1)$$
for any $X\in\Gamma(TM)$.

\parindent=8mm Recently, Cabrerizo et. al [3] extended the above result into a characterization for a slant submanifold in a contact metric manifold. In fact, they have obtained the following theorem.\\

\noindent
{\bf{Theorem 3.1 [3]}} {\it Let $M$ be a submanifold of an almost contact metric manifold $\bar M$ such that $\xi\in TM$. Then $M$ is slant if and only if there exists a constant $\lambda\in[0,1]$ such that}
$$T^2=\lambda(-I+\eta\otimes\xi).\eqno(3.2)$$
{\it Furthermore, in such a case, if $\theta$ is slant angle, then it satisfies that $\lambda=\cos^2\theta$.}

\parindent=8mm
Hence, for a slant submanifold $M$ of an almost contact metric manifold $\bar M$, the following relations are consequences of the above theorem.
$$g(TX,TY)=\cos^2\theta[g(X,Y)-\eta(X)\eta(Y)]\eqno(3.3)$$
$$g(FX, FY) =\sin^2\theta[g(X, Y)-\eta(X)\eta(Y)]\eqno(3.4)$$
for any $X, Y\in\Gamma(TM)$.

\parindent=8mm In the following theorem we consider $M$ as a totally umbilical slant submanifold of a Kenmotsu manifold $\bar M$.\\

\noindent
{\bf{Theorem 3.2}} {\it{Let $M$ be a totally umbilical slant submanifold of a Kenmotsu manifold $\bar M$. Then at least one of the following statements is true}}
\begin{enumerate}
\item [(i)] {\it{$M$ is invariant}}
\item [(ii)] {\it{$M$ is anti-invariant}}
\item [(iii)] {\it{$M$ is totally geodesic}}
\item [(iv)] $dim M=1$
\item [(v)] {\it{If $M$ is proper slant, then $H\in\Gamma(\mu)$}}
\end {enumerate}
{\it{where $H$ is the mean curvature vector of $M$.}}\\

\noindent{\it{Proof.}} As $M$ is totally umbilical slant submanifold, then we have
$$h(TX, TX)=g(TX, TX)H=\cos^2\theta\{\|X\|^2-\eta^2(X)\}H.$$
Using (2.5), we obtain
$$\cos^2\theta\{\|X\|^2-\eta^2(X)\}H=\bar\nabla_{TX}TX-\nabla_{TX}TX.$$
Then from (2.8), we get
$$\cos^2\theta\{\|X\|^2-\eta^2(X)\}H=\bar\nabla_{TX}\phi X-\bar\nabla_{TX}FX-\nabla_{TX}TX.$$
By (2.6) and (2.11), we derive
$$\cos^2\theta\{\|X\|^2-\eta^2(X)\}H=(\bar\nabla_{TX}\phi)X+\phi\bar\nabla_{TX}X+A_{FX}TX$$
$$~~~~~~~~~~~~~~~~~-\nabla^\perp_{TX}FX-\nabla_{TX}TX.$$
Using (2.3) and (2.5), we obtain
$$\cos^2\theta\{\|X\|^2-\eta^2(X)\}H=g(\phi TX, X)\xi-\eta(X)\phi TX+\phi(\nabla_{TX}X+h(X, TX))$$
$$~~~~~~~~~~~~~~~~~+A_{FX}TX-\nabla^\perp_{TX}FX-\nabla_{TX}TX.$$
From (2.8), (2.10), (2.17) and the fact that $X$ and $TX$ are orthogonal vector fields on $M$, we arrive at
$$\cos^2\theta\{\|X\|^2-\eta^2(X)\}H=-g(TX, TX)\xi-\eta(X)T^2X-\eta(X)FTX+T\nabla_{TX}X$$
$$~~~~~~~~~~~~~~~~~~~~~~~~~+F\nabla_{TX}X+A_{FX}TX-\nabla^\perp_{TX}FX-\nabla_{TX}TX.$$
Then, using (3.2) and (3.3), we get
$$\cos^2\theta\{\|X\|^2-\eta^2(X)\}H=-\cos^2\theta\{\|X\|^2-\eta^2(X)\}\xi-\cos^2\theta\eta(X)\{-X+\eta(X)\xi\}$$
$$~~~~~~~~~~-\eta(X)FTX+T\nabla_{TX}X+F\nabla_{TX}X$$
$$~~~~~~~~~~~~~~~~~+A_{FX}TX-\nabla^\perp_{TX}FX-\nabla_{TX}TX.\eqno(3.5)$$
Taking the inner product with $TX$ in (3.5), for any $X\in\Gamma(TM)$, we obtain
$$0=g(T\nabla_{TX}X, TX)+g(A_{FX}TX, TX)-g(\nabla_{TX}TX, TX).\eqno(3.6)$$
Now, we compute the first and last term of (3.6) as follows
$$g(T\nabla_{TX}X, TX)=\cos^2\theta\{g(\nabla_{TX}X, X)-\eta(X)g(\nabla_{TX}X,\xi)\}.\eqno(3.7)$$
Also, we have
$$g(\nabla_{TX}TX, TX)=g(\bar\nabla_{TX}TX, TX).$$
Using the property of Riemannian connection the above equation will be
$$g(\nabla_{TX}TX, TX)=\frac{1}{2}TXg(TX, TX)=\frac{1}{2}TX\{\cos^2\theta(g(X, X)-\eta(X)\eta(X))\}.$$
Again by the property of Riemannian connection, we derive
$$g(\nabla_{TX}TX, TX)=\cos^2\theta\{g(\bar\nabla_{TX}X, X)-\eta(X)g(\bar\nabla_{TX}X, \xi)\}$$
$$~-\cos^2\theta\eta(X)g(\bar\nabla_{TX}\xi, X).\eqno(3.8)$$
Using (2.4) and the fact that $X$ and $TX$ are orthogonal vector fields on $M$, the last term of (3.8) is identically zero, then by (2.5), we obtain
$$g(\nabla_{TX}TX, TX)=\cos^2\theta\{g(\nabla_{TX}X, X)-\eta(X)g(\nabla_{TX}X, \xi)\}.\eqno(3.9)$$
Thus, from (3.7) and (3.9), we get
$$g(T\nabla_{TX}X, TX)=g(\nabla_{TX}TX, TX).\eqno(3.10)$$
Using this fact in (3.6), we obtain
$$0=g(A_{FX}X, TX)=g(h(TX, TX), FX).$$
As $M$ is totally umbilical slant, then from (2.17) and (3.3), we get
$$0=\cos^2\theta\{\|X\|^2-\eta^2(X)\}g(H, FX).\eqno(3.11)$$
Thus, from (3.11), we conclude that either $\theta=\pi/2$, that is $M$ is anti-invariant  which part (ii) or the vector field $X$ is parallel to the structure vector field $\xi$, i.e., $M$ is $1-$dimensional submanifold which is fourth part of the theorem or $H\perp FX$, for all $X\in\Gamma(TM)$, i.e., $H\in\Gamma(\mu)$ which is the last part of the thorem or $H=0$, i.e., $M$ is totally geodesic which is (iii) or $FX=0,~\forall X\in\Gamma(TM)$, i.e., $M$ is invariant which is part (i). This proves the theorem completely.$~\blacksquare$

\parindent=8mm Now, if we consider $M$, a proper slant submanifold of a Kenmotsu manifold $\bar M$, then neither $M$ is invariant nor anti-invariant (by definition of proper slant) and also neither $dim M=1$. Hence, by the above result, only possibility is that $H\in\Gamma(\mu)$ for a totally umbilical proper slant submanifold. Thus, we prove the following main result.\\

\noindent
{\bf{Theorem 3.3}} {\it Every totally umbilical proper slant submanifold of a Kenmotsu manifold is totally geodesic.}\\

\noindent
{\it{Proof.}} Let $M$ be a totally umbilical proper slant submanifold of a Kenmotsu manifold $\bar M$, then for any $X, Y\in\Gamma(TM)$, we have
$$\bar\nabla_X\phi Y-\phi\bar\nabla_XY=g(\phi X, Y)\xi-\eta(Y)\phi X.$$
From (2.5) and (2.8), we obtain
$$\bar\nabla_XTY+\bar\nabla_XFY-\phi(\nabla_XY+h(X, Y))=g(TX, Y)\xi-\eta(Y)TX-\eta(Y)FX.$$
Again using (2.5), (2.6) and (2.8), we get
$$g(TX, Y)\xi-\eta(Y)TX-\eta(Y)FX=\nabla_XTY+h(X, TY)-A_{FY}X~~~~~~~~~~~~~~~~~~~~~~~~$$
$$~~~~~~~~~~~~~~~~~~~~~~~~~~~~~~~~~~~~~~~~+\nabla^\perp_XFY-T\nabla_XY-F\nabla_XY-\phi h(X, Y).$$
As $M$ is totally umbilical, then
$$g(TX, Y)\xi-\eta(Y)TX-\eta(Y)FX=\nabla_XTY+g(X, TY)H-A_{FY}X+\nabla^\perp_XFY$$
$$~~~~~~~~~~~~~~~~~~~~~~~~~~~~~~~~~~~~~~~~~-T\nabla_XY-F\nabla_XY-g(X, Y)\phi H.\eqno(2.12)$$
Taking the inner product with $\phi H$ in (3.12) and using the fact that $H\in\Gamma(\mu)$, we obtain
$$g(\nabla^\perp_XFY, \phi H)=g(X, Y)\|H\|^2.$$
Using (2.6) and the property of Riemannian connection, the above equation takes the form
$$g(FY, \nabla^\perp_X\phi H)=-g(X, Y)\|H\|^2.\eqno(3.13)$$
Now, for any $X\in\Gamma(TM)$, we have
$$\bar\nabla_X\phi H=(\bar\nabla_X\phi)H+\phi\bar\nabla_XH.$$
Using (2.3), (2.6), (2.8) and the fact that $H\in\Gamma(\mu)$, we obtain
$$-A_{\phi H}X+\nabla^\perp_X\phi H=-TA_HX-FA_HX+\phi\nabla^\perp_XH.\eqno(3.14)$$
Also, for any $X\in\Gamma(TM)$, we have
$$g(\nabla^\perp_XH, FX)=g(\bar\nabla_XH, FX)$$
$$~~~~~~~~~~~~~~~~~~~~~=-g(H, \bar\nabla_XFX).$$
Using (2.8), we get
$$g(\nabla^\perp_XH, FX)=-g(H, \bar\nabla_X\phi X)+g(H, \bar\nabla_XPX).$$
Then from (2.5) and (2.11), we derive
$$g(\nabla^\perp_XH, FX)=-g(H, (\bar\nabla_X\phi)X)-g(H, \phi\bar\nabla_XX)+g(H, h(X, PX)).$$
Using (2.3) and (2.17), the first and last term of right hand side of the above equation are identically zero and hence by (2.2), the second term gives
$$g(\nabla^\perp_XH, FX)=g(\phi H, \bar\nabla_XX).$$
Again, using (2.5) and (2.17), finally we obtain 
$$g(\nabla^\perp_XH, FX)=g(\phi H, H)\|X\|^2=0.$$
This means that 
$$\nabla^\perp_XH\in\Gamma(\mu).\eqno(3.15)$$
Now, taking the inner product in (3.14) with $FY$, for any $Y\in\Gamma(TM)$, we get
$$g(\nabla^\perp_X\phi H, FY)=-g(FA_HX, FY)+g(\phi\nabla^\perp_XH, FY).$$
Using (3.15), the last term of the right hand side of the above equation will be zero and then from (3.4), (3.13), we obtain
$$g(X, Y)\|H\|^2=\sin^2\theta\{g(A_HX, Y)-\eta(Y)g(A_HX, \xi)\}.\eqno(3.16)$$
Hence, by (2.7) and (2.17), the above equation reduces to
$$g(X, Y)\|H\|^2=\sin^2\theta\{g(X, Y)\|H\|^2-\eta(Y)g(h(X, \xi), H)\}.\eqno(3.17)$$
Since, for a Kenmotsu manifold $\bar M$,~$h(X, \xi)=0$, for any $X$ tangent to $\bar M$, thus we obtain
$$g(X, Y)\|H\|^2=\sin^2\theta g(X, Y)\|H\|^2.$$
Therefore, the above equation can be written as
$$\cos^2\theta g(X, Y)\|H\|^2=0.\eqno(3.18)$$
Since, $M$ is proper slant, thus from (3.18), we conclude that $H=0$ i.e., $M$ is totally geodesic in $\bar M$. This completes the proof of the theorem.$~\blacksquare$

\parindent=8mm Now, we give the following counter example of totally geodesic submanifold of $R^5$.\\

{\bf{Example 3.1}} Consider a $3-$dimensional proper slant submanifold with the slant angle $\theta\in[0, \pi/2]$ of $R^5$ with its usual Kenmotsu structure
$$x(u, v, t)=2(u\cos\theta, u\sin\theta, v, 0, t).$$
If we denote by $M$ a slant submanifold, then its tangent space $TM$ span by the vectors
$$e_1=\frac{\partial}{\partial u}+2\cos\theta(\frac{\partial}{\partial x^1}+y^1\frac{\partial}{\partial t})+2\sin\theta(\frac{\partial}{\partial x^2}+y^2\frac{\partial}{\partial t}),$$
$$e_2=\frac{\partial}{\partial v}=2\frac{\partial}{\partial y^1},~~~~~~e_3=\frac{\partial}{\partial t}=\xi.$$
Moreover, the vector fields
$$e^\star_1=-2\sin\theta(\frac{\partial}{\partial x^1}+y^1\frac{\partial}{\partial t})+2\cos\theta(\frac{\partial}{\partial x^2}+y^2\frac{\partial}{\partial t}),$$
$$e^\star_2==2\frac{\partial}{\partial y^2}$$
form the basis of $T^\perp M$. Furthermore, using Koszul's formula, we get $\bar\nabla_{e_i}e_i=-e_3=-\xi,~i=1, 2$ and when $i\neq j$, then $\bar\nabla_{e_i}e_j=0$, for  $i,j=1, 2, 3$. Also, $\bar\nabla_{e_3}e_3=0$, thus, from Gauss formula and (2.14), we obtain
$$h(e_1, e_1)=0,~~~~h(e_2, e_2)=0,~~~~h(e_1, e_2)=0,~~~h(e_3, e_3)=0$$
and hence we conclude that $M$ is totally geodesic.

Author' address:\\

Siraj UDDIN

\noindent Institute of Mathematical Sciences, Faculty of Science,University of Malaya, 50603 Kuala Lumpur, MALAYSIA

\noindent {\it E-mail}: {\tt siraj.ch@gmail.com}\\
\smallskip

Zafar AHSAN

\noindent Department of Mathematics, Aligarh Muslim University, 202002 Aligarh, INDIA

\noindent {\it E-mail}: {\tt zafar.ahsan@rediffmail.com}\\
\smallskip

Abdul Hadi YAAKUB

\noindent Institute of Mathematical Sciences, Faculty of Science,University of Malaya, 50603 Kuala Lumpur, MALAYSIA

\noindent {\it E-mail}: {\tt abdhady@um.edu.my}\\
\smallskip

\end{document}